\documentclass[12pt,reqno]{amsart}

\usepackage{graphicx}
     \makeatletter
     \def\section{\@startsection{section}{1}%
      \z@{.7\linespacing\@plus\linespacing}{.5\linespacing}%
     {\bfseries%\normalfont\scshape
     \centering
     }}
     \def\@secnumfont{\bfseries}
     \makeatother
\usepackage{amsthm}
\newtheorem*{theorem}{Theorem}
\newtheorem*{proposition}{Proposition}
\def \min   {\text {\rm min}}

\def \lim   {\text {\rm lim}}

\def \qed   {\hfill \vrule height6pt width 6pt depth 0pt}

\begin{document}
\title[]{Positivity properties of the matrix $\left[(i+j)^{i+j}\right]$}

\author[Rajendra Bhatia]{Rajendra Bhatia}
\address{Indian Statistical Institute\\ New Delhi 110016, India}
\address{Sungkyunkwan University, Suwon 440-746, Korea}

\email{rbh@isid.ac.in}

\author[Tanvi Jain]{Tanvi Jain}
\address{Indian Statistical Institute\\ New Delhi 110016, India}
\email{tanvi@isid.ac.in}

\subjclass[2010]{ 15B48, 42A82}

\keywords{Positive definite, totally positive, infinitely divisible, conditionally positive definite.}

% \begin{left}
% (Work in progress)\\
% May 2014
% \end{left}

\begin{abstract}
Let $p_1<p_2<\cdots<p_n$ be positive real numbers. It is shown that the matrix whose $i,j$ entry is $(p_i+p_j)^{p_i+p_j}$ is infinitely divisible, nonsingular and totally positive.
\end{abstract}

\maketitle

{\bf{1. Introduction}}\\

Matrices whose entries are obtained by assembling natural numbers in special ways often possess interesting properties. The most famous example of such a matrix is the Hilbert matrix $ H=\left[\frac{1}{i+j-1}\right]$ which has inspired a lot of work in diverse areas. Some others are the min matrix $M=\begin{bmatrix}\min (i,j)\end{bmatrix},$ and the Pascal matrix  $P=\left[\binom{i+j}{i}\right]$. There is a considerable body of literature around each of these matrices, a sample of which can be found in \cite{RBHi}, \cite{rb4} and \cite{MDC}.\\

In this note we initiate the study of one more matrix of this type. Let $A$ be the $ n \times  n$ matrix with its $(i,j)$ entry equal to $(i+j-1)^{i+j-1}$. Thus
\begin{equation}
A=\left[
\begin{array}{ccccc}
1 & 2^2 & 3^3 &\cdots & n^n\\
2^2 & 3^3 & 4^4 &\cdots &( n+1)^{n+1}\\
3^3 & 4^4 & 5^5 &\cdots &\cdots\\
\cdots &\cdots &\cdots &\cdots &\cdots\\
n^n &\cdots &\cdots &\cdots & (2n-1)^{2n-1}\\
\end{array}\right]. \label{eq1}
\end{equation}
More generally, let $p_1<p_2<\cdots <p_n$ be positive real numbers, and consider the $ n \times  n$  matrix
\begin{equation}B=\left[(p_i+p_j)^{p_i+p_j}\right]. \label{eq2}
\end{equation}
The special choice $ p_i=i-1/2$ in \eqref{eq2} gives us the matrix \eqref{eq1}. We investigate the behaviour of these matrices with respect to different kinds of positivity.\\

A real symmetric matrix $S$ is said to be {\it positive semidefinite} (psd) if  for every vector 
$x$, we have $\langle x, Sx \rangle \ge 0$. Further if  $\langle x, Sx\rangle =0$ only when $x = 0,$ then we say $S$ is {\it positive definite}. This is equivalent to saying that S is psd and nonsingular. If $S$ is a psd matrix, then for every positive integer $m,$ the $m$th {\it Hadamard power} (entrywise power) $S^{\circ m}=\left[s_{ij}^m\right]$  is also psd. Now suppose $ s_{ij}\ge0$. We say that $ S$ is {\it infinitely divisible} if for every real number $ r>0$, the matrix $S^{\circ r}=\left[s_{ij}^r\right]$   is psd. (See \cite{RBHi}, Chapter 5 of \cite{rbh}, and Chapter 7 of \cite{HJ} for expositions of this topic.) The principal minors of a psd matrix are nonnegative. This may not be so for other minors. A matrix with nonnegative entries is called {\it totally positive} if  all its minors are nonnegative. It is called {\it strictly totally positive} if all its minors are positive. We recommend the books \cite{fj, k, AP} and the survey article \cite{a} for an account of totally positive matrices.

\noindent Our main result is the following: 

\begin{theorem}
Let $p_1< p_2<\cdots <p_n$ be positive real numbers. Then the matrix $B$ defined in \eqref{eq2} is infinitely divisible, nonsingular and totally positive.
\end{theorem}

There is another way of stating this. Let $X$ be a subset of $\mathbb{R}$. A continuous function $K:X\times X\to\mathbb{R}$ is said to be a {\it positive definite kernel} if for every $n$ and for every choice $x_1<\cdots<x_n$ in $X$, the matrix $\begin{bmatrix} K(x_i,x_j)\end{bmatrix}$ is positive definite. In the same way we can define infinitely divisible and totally positive kernels. Our theorem says that the kernel $K(x,y)=(x+y)^{x+y}$ on $(0,\infty)\times (0,\infty)$ is infinitely divisible and totally positive. This is an addition to the examples given in \cite{RBHi, fj, k, AP}. The three matrices in the first paragraph also have the properties mentioned in the theorem.\\

{\bf 2. {Proof}}\\

Let $H_1$ be the space of all vectors $ x=(x_1,...,x_n)$ with $\sum x_i = 0$. A real symmetric matrix $S$ is said to be {\it conditionally positive definite} (cpd) if $\langle x,Sx\rangle \ge 0$ for all $ x\in  H_1$. If $-S$ is cpd, then $S$ is said to be {\it conditionally negative definite} (cnd). According to a theorem of C. Loewner, a matrix $S=[s_{ij}]$ is infinitely divisible if and only if the matrix $[\log  s_{ij}]$ is cpd. See Exercise 5.6.15 in \cite{rbh}.
\vskip0.2in

By Loewner's theorem cited above, in order to prove that the matrix $B$ defined in \eqref{eq2} is infinitely divisible it is enough to show that the matrix
\begin{equation}
C=\left[(p_i+p_j) \log(p_i+p_j)\right]
 \label{eq3}
\end{equation}
is cpd. It is convenient to use the formula
$$ \log x=\int_0^{\infty}\left(\frac{1}{1+\lambda}-\frac{1}{x+\lambda}\right)d \lambda,  \,\,\, x>0,$$
which can be easily verified. Using this we can write our matrix $C$ as
\begin{equation*}
 C=\left[\int_0^{\infty}\left(\frac{p_i+p_j}{1+\lambda}-\frac{p_i+p_j}{p_i+p_j+\lambda}\right)d \lambda\right]. 
 \end{equation*}
 We will show that the matrix $[p_i+p_j]$ is cpd, and the matrix $\left[\frac{p_i+p_j}{p_i+p_j+\lambda}\right]$ is cnd for each $\lambda>0$. From this it follows that $C$ is a cpd matrix.
\vskip0.2in
\noindent Let $ D$ be the diagonal matrix $ D=\mbox{diag}(p_1,...,p_n)$ and $E$ the matrix with all its entries equal to 1. Then $[p_i+p_j]=DE + ED.$
Every vector $x$ in $H_1$ is annihilated by $E.$ Hence $\langle x,(DE+ED)x \rangle=\langle x,DEx\rangle+ \langle Ex, Dx \rangle=0.$  So, the matrix $[p_i+p_j]$ is cpd. Using the identity
\begin{equation*}
\frac{p_i+p_j}{p_i+p_j+\lambda}=1-\frac{\lambda}{p_i+p_j+\lambda},
\end{equation*}
we can write
\begin{equation*}
\left[\frac{p_i+p_j}{p_i+p_j+\lambda}\right]=E-\lambda C_{\lambda},
\end{equation*}
where $C_{\lambda}=\left[\frac{1}{p_i+p_j+\lambda}\right]$. This is a Cauchy matrix (see \cite{rbh}) and is positive definite. Hence it is also cpd. The matrix $E$ annihilates $H_1$, and  therefore is cnd. Hence $E- \lambda C_\lambda$ is cnd for every $\lambda >0$. This completes the proof of the assertion that $C$ is cpd, and $B$ infinitely divisible.

\vskip0.2in
Since  $C_\lambda$ is positive definite, $\langle x,C_\lambda x \rangle>0$ for every non zero vector $x.$   If $x \in H_1$, then $Ex=0$, and $\langle x, (DE+ED)x\rangle =0.$ So, the arguments  given above also show that $\langle x, Cx \rangle >0 $ for every non zero vector $x$ in $H_1$. By  Lemma 4.3.5  in \cite{BR}, this condition is necessary and sufficient for $C$ to be nonsingular. Using the next proposition, we can conclude that $B$ is nonsingular.

\begin{proposition}
If $C$ is a nonsingular conditionally positive definite matrix, then the matrix $[e^{c_{ij}}]$ is positive definite.
\end{proposition}

\noindent {\it Proof.} By Proposition 5.6.13 of \cite{rbh} we can express C as 
\begin{equation*}
C= P+ YE + E\overline{Y},
\end{equation*}
where $P$ is a psd matrix and $Y$ is a diagonal matrix. By Problem 7.5.P.25 in \cite{HJ} the matrix $[e^{c_{ij}}]$ is  positive definite unless $P$ has two equal columns. Suppose the $i$th column of $P$ is equal to its $j$th column. Let 
$x$ be any vector with coordinates $ x_i=-x_j\neq 0$, and all other coordinates zero. Then $x \in H_1$ and $Px=0.$  Hence $\langle x, Cx \rangle=0.$ This is not possible since $C$ is a nonsingular cpd matrix.\hfill{$\qed$}

\vskip.2in
We have proved that the matrix $B$ is infinitely divisible and nonsingular. These properties are inherited by the matrix $A$ defined in \eqref{eq1}. This is, moreover, a {\it Hankel matrix}; i.e, each of its  antidiagonals has the same entry. Theorem 4.4 of \cite{AP} gives a simple criterion for strict total positivity of such a matrix.  According to this a Hankel matrix $A$ is strictly totally positive if  and only if $A$ is positive definite and so is the matrix $\tilde{A}$ obtained  from $A$  by deleting its first column and last row. For the matrix $A$ in (1), $\tilde{A}$ is the $(n-1) \times (n-1)$ matrix whose $(i,j)$ entry  is $(i+j)^{i+j}$. Both $A$ and $\tilde{A}$ are positive definite. Hence $A$ is strictly totally positive. In fact we have shown that for every $r>0$, the matrix $A^{\circ r}$ is strictly totally positive.\\

Now let $k_1<k_2<\cdots<k_n$ be positive integers. The matrix $K$ with entries $k_{ij}=(k_i+k_j)^{k_i+k_j}$ is principal submatrix of $A$. Hence it is infinitely divisible and strictly totally positive. The same holds for $K^{\circ r}$ for every $r>0$. Next let $0<q_1<q_2<\cdots<q_n$ be rational numbers. Let $q_j=l_j/m_j$, where $l_j$ and $m_j$ are positive integers. Let $m$ be the LCM of $m_1,\ldots, m_n$ and $k_j=mq_j$. Then $k_1<k_2<\cdots<k_n,$ and as seen above, the matrix $K=\begin{bmatrix}(k_i+k_j)^{k_i+k_j}\end{bmatrix}$ is infinitely divisible and strictly totally positive. Now consider the matrix $Q=\begin{bmatrix}(q_i+q_j)^{q_i+q_j}\end{bmatrix}.$ Then for each $r>0$
\begin{eqnarray*}
Q^{\circ r} & = & \begin{bmatrix}(q_i+q_j)^{(q_i+q_j)r}\end{bmatrix}\\
 & = & \begin{bmatrix}\frac{(k_i+k_j)^{(k_i+k_j)r/m}}{m^{q_i r}m^{q_j r}}\end{bmatrix}\\
 & = & XK^{\circ r/m}X^*,
\end{eqnarray*}
where $X$ is the positive diagonal matrix with entries $\frac{1}{m^{q_1r}}, \ldots, \frac{1}{m^{q_nr}}$ on its diagonal. We have seen that the matrix $K^{\circ r/m}$ is positive definite and strictly totally positive. Hence, so is the matrix $Q^{\circ r}.$ A continuity argument completes the proof of the theorem.\hfill{$\qed$}

\vskip.2in
We believe that the matrix $B$ in \eqref{eq2} is strictly totally positive. However the continuity argument that we have invoked at the last step only shows that it is a limit of such matrices.\\

Like for the other matrices mentioned in the opening paragraph, it would be interesting to have formulas for the determinant of $A$.\\
 
In a recent work \cite{rj} of ours, we have studied spectral properties of the matrices $\left[(p_{i}+p_{j})^r\right]$, where $r$ is any positive real number.

\vskip.2in
{\it We thank R. B. Bapat for illuminating discussions, and in particular for the argument in the last paragraph of our proof. The first author is supported by a J. C. Bose National Fellowship and was a Fellow Professor of Sungkyunkwan University in the summer of 2014. The second author is supported by a SERB Women Excellence Award.}

\end{document}